\documentclass[10pt,a4paper,reqno]{amsart}
\newtheorem{theorem}{\bf Theorem}[section]
\newtheorem{remark}[theorem]{Remark}
\newtheorem{prop}[theorem]{Proposition}
\newtheorem{conj}[theorem]{Conjecture}
\newtheorem{cor}[theorem]{Corollary}
\usepackage{xcolor}
\usepackage{amsmath}
\usepackage[colorlinks=true,linkcolor=blue,citecolor=orange]{hyperref}

\DeclareMathAlphabet{\mathpzc}{OT1}{pzc}{m}{it}

\begin{document}

\title[On discrete X-ray transform
]{
On 
discrete X-ray
transform}

\author{Roman Novikov}
\address{Roman G. Novikov, CMAP, CNRS, \'{E}cole polytechnique, Institut Polytechnique de Paris, 91128 Palaiseau, France \newline
\& IEPT RAS, 117997 Moscow, Russia}
\email{novikov@cmap.polytechnique.fr}

\author{Basant Lal Sharma}
\address{Basant Lal Sharma, Department of Mechanical Engineering, Indian Institute of Technology Kanpur, Kanpur, 208016 UP, India}
\email[Corresponding author]{bls@iitk.ac.in}

\date{\today}

\begin{abstract}
We consider a discrete version of X-ray 
transform
going back, in particular, to Strichartz (1982).
We suggest non-overdetermined reconstruction for this discrete transform.
Extensions to weighted (attenuated) analogues are given.
Connections to the continuous case are presented.

\bigskip



\end{abstract}

\maketitle 

\section{Introduction}

We consider 
the $X$-ray transform $\mathcal{P}$ along oriented straight lines $\gamma$ in $\mathbb{R}^d,~~ d\ge2,$ defined by 
\begin{equation}
\mathcal{P}f(\gamma)=\int_{\gamma} f(y)dy,\quad \gamma=(\theta,x)\in T,
\label{XTeq1}
\end{equation}
with $T$ defined by
\begin{equation}
T=\{(\theta,x)\in\mathbb{S}^{d-1}\times\mathbb{R}^d: x\cdot\theta=0\},
\label{XTeq2}
\end{equation}
where $f$ is a test function on $\mathbb{R}^d$.
Here,
$\gamma=(\theta,x)$ is considered as the oriented straight line (ray) in $\mathbb{R}^d$ defined by
\begin{equation}
\gamma=(\theta,x)=\{y \in \mathbb{R}^d: y=x+t\theta,\quad t\in\mathbb{R}\},
\label{XTeq3}
\end{equation}
where $\theta$ gives the orientation of $\gamma.$
In connection with the above definition
and formulas for finding $f$ from $\mathcal{P}f$, see, for example, \cite{Radon} and \cite{Natterer}.

Note that 
\begin{equation}
\text{dim }T=2d-2.
\end{equation}
Therefore, the problem of finding $f$ from $\mathcal{P}f$ on $T$ is not over-determined for $d=2$ and is over-determined for $d\ge3$.
Therefore, finding $f$ from $\mathcal{P}f$ on $\Gamma\subset T,$ where dim $\Gamma=d$ is of main interest for $d\ge3$.

In this work, we continue studies on discrete analogues of the X-ray transform $\mathcal{P}$ in \eqref{XTeq1}.
In connection with the relevant literature, see, for example, \cite{Strichartz}, \cite{Beylkin}, \cite{Stein}, \cite{Fokas}, \cite{Averbuch}, \cite{Press}, \cite{Kawazoe}.
More precisely, the present work is focused on discrete X-ray transforms
going back to Strichartz \cite{Strichartz} for $d=2$.

In this context, in addition to the transform $\mathcal{P}$ defined by \eqref{XTeq1}, we also consider
its discrete counterpart, 
for functions $f$ on $\mathbb{Z}^d\subset\mathbb{R}^d$,
defined by
\begin{equation}
\mathcal{P}f(\gamma)=\sum_{y\in\gamma\cap\mathbb{Z}^d} f(y),\quad \gamma\in {T}',
\label{XTeq1disc}
\end{equation}

\begin{equation}
{T}'=\cup_{\zeta\in \mathbb{Z}^d} {T}_{\zeta},\quad {T}_{\zeta}:=\{\gamma\in T: \zeta\in\gamma\}.
\label{defTTz}
\end{equation}

One can see that ${T}'$ is not a countable set in contrast to $\mathbb{Z}^d$.
Moreover, we have already that
\begin{equation}
{T}_\zeta\simeq \mathbb{S}^{d-1},\quad
\zeta\in \mathbb{Z}^d.
\end{equation}

Therefore, in the discrete case \eqref{XTeq1disc}, it is natural to consider $\mathcal{P}f$ restricted to countable subsets
$\Gamma\subset T'.$

In the present work, in the discrete case, our basic assumption on $f$ is that
\begin{equation}
\text{supp }f\subseteq B_r\cap \mathbb{Z}^d,
\label{suppfBr}
\end{equation}
where
\begin{equation}
B_r:=\{x\in\mathbb{R}^d: |x|\le r\}, \quad r>0.
\label{diskBr}
\end{equation}
One can see that
\begin{equation}
N_r:=\#(B_r\cap\mathbb{Z}^d)=\Theta(r^{d}),\quad r\to\infty;
\label{defNrn}
\end{equation}
or equivalently, $c\, r^{d}\le N_r\le c'\,r^{d}$ for some $c,c'>0$, as $r\to\infty$.
Moreover, 
$N_r=|B_r|+O(r^{d-1})$, as $r\to\infty$, where $|.|$ denotes the volume; this result, together with more precise estimates for the remainder, goes back to Gauss, at least, for $d=2$.
See, for example, \cite{Kratzel}.

Therefore, 
in the discrete case \eqref{XTeq1disc}, \eqref{suppfBr}, it is especially interesting
to consider the problem of finding $f$ from $\mathcal{P}f$ restricted to $\Gamma$ such that 
$\#\Gamma=N_r$ or, at least, $\#\Gamma=\Theta(r^d)$.
In this work, we are focused on such non-overdetermined reconstructions.

Let
\begin{equation}
\Omega=\{\theta: \theta=z/|z|\text{ for some } z\in\mathbb{Z}^d\setminus\{0\}\}.
\label{defOmg}
\end{equation}
We say that a ray $\gamma\in T'$ is {\em rational} if $\hat{\gamma}\in\Omega$
and we say that
$\gamma\in T'$ is {\em irrational} if $\hat{\gamma}\notin\Omega$, where $\hat{\gamma}$ denotes the direction of $\gamma$.

In some sense, the simplest way to suggest a non-overdetermined inversion of the transform $\mathcal{P}$ in \eqref{XTeq1disc}
consists in considering 
\begin{equation}
\Gamma=\Gamma_\theta:=\{\gamma\in {T}': \hat{\gamma}=\theta\}, \quad\theta\in \mathbb{S}^{d-1}\setminus\Omega.
\label{newdef}
\end{equation}
In this case, the non-overdetermined reconstruction of $f$ from $\mathcal{P}f$ on $\Gamma_\theta$ is given by the obvious formula \eqref{fdef1} in subsection \ref{sec2p1} of the present work.

On the other hand, in the discrete case \eqref{XTeq1disc}, it is very natural to consider the problem of finding $f$ from $\mathcal{P}f$ restricted to the set $\mathpzc{T}$ of rational rays, where 
\begin{equation}
\mathpzc{T}:=\{\gamma\in {T}': \hat{\gamma}\in\Omega\}.
\label{defTset}
\end{equation}
The definition of $\mathcal{P}f$ on $\mathpzc{T}$ goes back to Strichartz \cite{Strichartz}, at least, for $d=2$.
However, the problem of finding $f$ on $\mathbb{Z}^d$, from $\mathcal{P}f$ on $\mathpzc{T}$, e.g., under assumption \eqref{suppfBr}, is 
overdetermined in the sense that 
\begin{equation}
\# \mathpzc{T}_r^{\max}=\infty,
\label{defTinf}
\end{equation}
\begin{equation}
\# \mathpzc{T}_r^{\min}=\Theta((r^{d})^2)\text{ as } r\to\infty,
\label{defT}
\end{equation}
in comparison with $N_r$ in \eqref{defNrn}, where
\begin{equation}
\mathpzc{T}_r^{\max}:=\{\gamma\in \mathpzc{T}: \gamma\cap B_r\cap\mathbb{Z}^d\ne\emptyset\},
\label{defTrmax}
\end{equation}
\begin{equation}
\mathpzc{T}_r^{\min}:=\{\gamma\in \mathpzc{T}: \#(\gamma\cap B_r\cap\mathbb{Z}^d)\ge2\}.
\label{defTrmin}
\end{equation}

For completeness of presentation, we prove \eqref{defT} in subsection \ref{secproof4p1}.
The lower bound in \eqref{defT}
is not obvious 
and the proof 
in
\ref{secproof4p1} 
uses, in particular, estimates for the Farey points.

In the present work, for the case of rational rays, 
we provide a non-overdetermined reconstruction of Cormack-type of $f$ from $\mathcal{P}f$ on an appropriate $\Gamma\subset\mathpzc{T}$, under assumption \eqref{suppfBr}; see Theorem \ref{thm2p1} and Theorem \ref{thm2p2} in subsection \ref{secrat}.
In fact, some non-overdetermined reconstruction for this case was already suggested 
in \cite{Kawazoe} for $d=2$ but without clarifying its non-overdetermined property.
In subsection \ref{secrat}, in connection with Theorem \ref{thm2p1},
we also formulate Conjectures \ref{conjec2p31} and \ref{conjec2p3}, as discrete analogues of well known results in the continuous case, where $f$ is not compactly supported.

The aforementioned Theorems \ref{thm2p1} and \ref{thm2p2} 
admit straightforward generalizations to the 
weighted (attenuated) case; see subsection \ref{secwt2p1}.

Note that the discrete X-ray transform $\mathcal{P}$ in \eqref{XTeq1disc}
arises, in particular, in the framework of continuous one in \eqref{XTeq1} for the case when 
\begin{equation}
f(x)=\sum_{y\in\mathbb{Z}^d}c(y)\delta(x-y),
\label{fexpand}
\end{equation}
where $\delta$ denotes the Dirac-delta function; see subsection \ref{sec3p1}.
Besides this, the discrete X-ray transform $\mathcal{P}$ in \eqref{XTeq1disc} can be applied to studies of
the continuous one in \eqref{XTeq1} for the case of regular $f$ as sketched in subsection \ref{sec3p2}.

\section{Non-overdetermined inversions for discrete X-ray transform}
\label{sec2}

\subsection{Irrational rays}
\label{sec2p1}
The simplest way to consider non-overdetermined inversion of the discrete X-ray transform $\mathcal{P}$ in \eqref{XTeq1disc} consists in the following formula for finding $f$ on $\mathbb{Z}^d$ from $\mathcal{P}f$ on $\Gamma_\theta$ as in \eqref{newdef}:
\begin{equation}
f(x)=\mathcal{P}f(\gamma_{x,\theta}), \quad x\in\mathbb{Z}^d,
\label{fdef1}
\end{equation}
where ${\gamma}_{x,\theta}$ denotes the straight line ${\gamma}$ such that $x\in\gamma, \hat{\gamma}=\theta.$
This formula follows from 
\eqref{XTeq1disc} and the observation that, for fixed $x\in\mathbb{Z}^d$ and $\theta\in\mathbb{S}^{d-1}\setminus\Omega$, we have
$\gamma_{x,\theta}\cap\mathbb{Z}^d=x.$
One can see also
that $\#(\gamma\cap\mathbb{Z}^d)=1$ if $\gamma\in\Gamma_\theta$
and
that there is one-to-one correspondence between $\Gamma_\theta$ in \eqref{newdef} and $\mathbb{Z}^d$.

Formula \eqref{fdef1} admits the following generalization:
\begin{equation}
f(x)=\mathcal{P}f(\gamma_{x,\theta(x)}), \quad x\in\mathbb{Z}^d,\quad \theta(x)\in\mathbb{S}^{d-1}\setminus\Omega,
\label{fdeff2}
\end{equation}
where $\gamma_{x,\theta(x)}$ is defined as in \eqref{fdef1} but with $\theta=\theta(x)$ depending on $x$.

Note that formulas \eqref{fdef1}, \eqref{fdeff2} do not require condition \eqref{suppfBr}
and, in principle, are fulfilled by any complex valued function $f$ on $\mathbb{Z}^d$.

Let
\begin{equation}
\Omega_\rho=
\{\theta: \theta=z/|z|,\,
z=(z_1, \dotsc, z_d)\in\mathbb{Z}^d\setminus\{0\},\,
|z|/\text{gcd}(z_1, \dotsc, z_d)\le \rho\},
\label{omgrho}
\end{equation}
where
$\rho\ge1$, and gcd stands for the greatest common divisor.

In addition to formula \eqref{fdeff2}, under condition \eqref{suppfBr}, the following formula also holds:
\begin{equation}
f(x)=\mathcal{P}f(\gamma_{x,\theta(x)}), \quad x\in\mathbb{Z}^d,\, \theta(x)\in\mathbb{S}^{d-1}\setminus\Omega_\rho,\, \rho>2r.
\label{fdeff2n}
\end{equation}

Note that
\eqref{fdeff2n} permits all irrational rays $\gamma_{x,\theta(x)}$ and the rational rays $\gamma_{x,\theta(x)}$ with a sufficiently high degree of irrationality.
However, it is specially natural and interesting to consider the discrete X-ray transform $\mathcal{P}$ in \eqref{XTeq1disc} for rational rays $\gamma\in\mathpzc{T}$ only, in a way similar to the discrete Radon transform suggested in \cite{Strichartz}.
Moreover, under 
condition \eqref{suppfBr}, it is natural to study reconstruction of $f$ from $\mathcal{P}f$
on $\Omega_\rho$ with $\rho=r$.

\subsection{Rational rays}
\label{secrat}
The simplest way to construct a subset $\mathpzc{T}^\star\subset\mathpzc{T}$,
where $\mathpzc{T}$ is defined by \eqref{defTset},
and $\mathpzc{T}^\star$ is in one-to-one correspondence with $\mathbb{Z}^d$,
is as follows.

Let
\begin{equation}
\gamma_z=(\hat{z}_\star,z)=\{y \in \mathbb{R}^d: y=z+t \hat{z}_\star,\quad t\in\mathbb{R}\}, \quad z\in \mathbb{Z}^d, \quad z_1^2+z_2^2\ne0,
\label{RTeqHz1}
\end{equation}
\begin{equation}
\hat{z}_\star=\frac{-z_2 e_1+z_1e_2}{\sqrt{z_1^2+z_2^2}}, 
\quad z=\sum_{i=1}^dz_ie_i,
\label{RTeqHz1s}
\end{equation}
and
\begin{equation}
\gamma_z=(\theta,z), \quad \theta\in\Omega, \quad\theta=(\theta_1,\theta_2,0,\dotsc,0),\quad z\in \mathbb{Z}^{d},\quad z_1^2+z_2^2=0.
\label{RTeqHz2}
\end{equation}
Here, $e_1, \dotsc, e_d$ are the standard basis vectors in $\mathbb{R}^d.$

One can assume that $\theta$ is fixed in \eqref{RTeqHz2}.
Let
\begin{equation}
\mathpzc{T}^{\star}=\cup_{z\in \mathbb{Z}^d}\gamma_z.
\label{defTstar}
\end{equation}
We also consider
\begin{equation}
\mathpzc{A}_{\alpha,\beta}=\{z\in\mathbb{Z}^d: \alpha \le \sqrt{z_1^2+z_2^2}\le \beta\},
\label{defCab}
\end{equation}
\begin{equation}
\mathpzc{T}^{\star}_{\alpha,\beta}=\{\gamma_z\in\mathpzc{T}^{\star}: \alpha\le \sqrt{z_1^2+z_2^2}\le \beta\},\quad 0\le \alpha\le \beta,
\label{defTstarab}
\end{equation}
where
$z=(z_1, z_2, \dotsc, z_d).$

\begin{theorem}
 Let $f$ satisfy \eqref{suppfBr}.
Then $\mathcal{P}f$ on $\mathpzc{T}^\ast$ uniquely determines $f$,
where $\mathpzc{T}^\ast$ is defined by \eqref{defTstar}.
Moreover, this determination is given via the two-dimensional formulas
\eqref{defrecon1}-\eqref{defrecon2}
for $d=2$,
and adopted for each slice with fixed $z_3, \dotsc, z_d$ for $d\ge 3$.
In addition, in this framework,
$\mathcal{P}f$ on $\mathpzc{T}^{\star}_{\alpha,\beta}$,
where $\beta\ge r$, uniquely determines $f$ on $\mathpzc{A}_{\alpha,\beta}$.
\label{thm2p1}
\end{theorem}

\subsection*{Proof of Theorem \ref{thm2p1}}
\label{secthm2p1}
We start with the case $d=2.$

Consider the sets
$S_{j}$,
$j=1, 2, \dotsc J,$
defined as follows:
\begin{equation}
S_1=\{z\in B_r\cap \mathbb{Z}^2: |z|\ge |\zeta|\text{ for any }\zeta\in B_r\cap \mathbb{Z}^2\},\quad
S_J=\{0\},
\label{defXS1}
\end{equation}
\begin{equation}
S_{i+1}=\{z\in (B_r\cap \mathbb{Z}^2)\setminus \cup_{k=1}^i S_k: |z|\ge |\zeta|\text{ for any }\zeta\in (B_r\cap \mathbb{Z}^2)\setminus \cup_{k=1}^i S_k\},
\label{defXSj}
\end{equation}
where $i+1\le J$, and $B_r$ is the disc in \eqref{suppfBr},\eqref{diskBr}.
One can see that
\begin{equation}
B_r\cap\mathbb{Z}^2=\cup_{j=1,\dotsc, J} S_j,\quad S_i \cap S_j=\emptyset\text{ if }i\ne j.
\label{Brcov}
\end{equation}

We reconstruct $f$ on $S_1$ first.
Then we reconstruct $f$ on $S_{i+1}$ inductively.
The reconstruction formulas are as follows:
\begin{equation}
f(z)=\mathcal{P}f(\gamma_z), \quad z\in S_1,
\label{defrecon1}
\end{equation}
\begin{equation}
f(z)=\mathcal{P}f(\gamma_z)-\sum_{\zeta\in \gamma_z\cap (\cup_{k=1}^i S_k)}f(\zeta), \quad z\in S_{i+1},
\label{defrecon2}
\end{equation}
for $i=1, 2, \dotsc, J-1,$
where
$\gamma_z$ is defined by \eqref{RTeqHz1}-\eqref{RTeqHz2}.
In view of \eqref{Brcov}, formulas \eqref{defrecon1}, \eqref{defrecon2} give a reconstruction on the whole $B_r\cap\mathbb{Z}^2.$

Formulas \eqref{defrecon1}, \eqref{defrecon2} follow from the property
that $z$ is the point of $\gamma_z$ closest to the origin $\{0\}$, definition of $\mathcal{P}f$ in \eqref{XTeq1disc},
and definitions of $S_j$ in \eqref{defXS1}, \eqref{defXSj}.
In particular, one can see that if $z\in S_i$, $z'\in \gamma_z$, $z'\ne z$, and $z'\in B_r\cap \mathbb{Z}^2$,
then $z'\in S_{j},\, j<i.$

In addition, formulas \eqref{defrecon1}, \eqref{defrecon2} give a reconstruction of $f$ on $\mathpzc{A}_{\alpha,\beta}$ from $\mathcal{P}f$ on $\mathpzc{T}^\star_{\alpha,\beta}$, $r\le \beta$.

This completes the proof for $d=2.$

In dimension $d\ge3,$ we proceed in a similar way for each two dimensional plane $\Xi=\{x\in\mathbb{R}^d: x=(x_1, x_2, z_3, \dotsc, z_d), (x_1, x_2)\in\mathbb{R}^2\}, \, z_3,\dotsc, z_d\in\mathbb{Z}$, where the origin $O$ (of $\Xi$) is considered at the point $x_0=(0, 0, z_3, \dotsc, z_d)$.
In this setting, our considerations on $\Xi$ reduce to our considerations for $d=2.$
This completes the proof of Theorem \ref{thm2p1}.

 One can see that the straight lines $\gamma\in \mathpzc{T}^\star$ used in Theorem \ref{thm2p1} are parallel to the plane $\Xi_{e_1,e_2}=$ span $\{e_1, e_2\}\text{ in } \mathbb{R}^d$. Theorem \ref{thm2p1} admits a straightforward generalization to the case of rational straight lines parallel to the plane 
\begin{equation}
\Xi_{a,b}=\text{ span }\{a, b\}\text{ in } \mathbb{R}^d,
\label{defXi}
\end{equation}
where $a, b\in \mathbb{Z}^d$, and $a, b$ are linearly independent. In this case, we define $\gamma_z$ by the formulas
\begin{equation}
\gamma_z=(\hat{z}_\star,z)=\{y \in \mathbb{R}^d: y=z+t \hat{z}_\star,\quad t\in\mathbb{R}\}, \quad z\in \mathbb{Z}^d,\quad
(z\cdot a)^2+(z\cdot b)^2\ne0,
\label{RTeqHz1new}
\end{equation}
\begin{equation}
\hat{z}_\star=\frac{-(z\cdot b)a+(z\cdot a) b}{|-(z\cdot b)a+(z\cdot a) b|}, 
\label{RTeqHz1snew}
\end{equation}
and
\begin{equation}
\gamma_z=(\theta,z), \quad \theta\in\Omega\cap\Xi_{a,b}, \quad z\in \mathbb{Z}^d, \quad (z\cdot a)^2+(z\cdot b)^2=0.
\label{RTeqHz2new}
\end{equation}
One can assume that $\theta$ is fixed in \eqref{RTeqHz2new}.

We also consider
\begin{equation}
\mathbb{Z}^2_{a,b}:=\mathbb{Z}^d\cap\Xi_{a,b},
\label{defXi2}
\end{equation}
and
\begin{equation}
\mathpzc{A}_{\alpha,\beta, a,b}=\{z\in\mathbb{Z}^d: \alpha \le |z_{a,b}|\le \beta\},
\label{defCab1}
\end{equation}
\begin{equation}
\mathpzc{T}^{\star}_{\alpha,\beta,a,b}=\{\gamma_z\in\mathpzc{T}^{\star}: \alpha\le|z_{a,b}|\le \beta\},\quad 0\le \alpha\le \beta,
\label{defTstarab1}
\end{equation}
where
\begin{equation}
z=z_{a,b}+z^{a,b}, \quad z_{a,b}, z^{a,b}\in\mathbb{R}^d,\quad z_{a,b}\in \Xi_{a,b}, \quad z_{a,b}\cdot z^{a,b}=0.
\end{equation}

\begin{theorem}
Let $f$ satisfy \eqref{suppfBr}
and $\mathpzc{T}^\ast$ be defined by \eqref{defTstar}, \eqref{RTeqHz1new}-\eqref{RTeqHz2new}.
Then $\mathcal{P}f$ on $\mathpzc{T}^\ast$ uniquely determines $f$.
Moreover, this determination is given via 
formulas \eqref{defrecon1}-\eqref{defrecon2}
for $d=2$,
where $\gamma_z$ are defined by \eqref{RTeqHz1new}-\eqref{RTeqHz2new},
and $S_j$ are defined by \eqref{defXS1}, \eqref{defXSj} with $\mathbb{Z}^2$ replaced by $\mathbb{Z}^2_{a,b}$,
and via these formulas adopted, for $d\ge 3$, for each two dimensional slice $\Xi$ in $\mathbb{R}^d$ such that $\Xi=\Xi_{a,b}+x_0$, where $x_0\in\mathbb{R}^d,\, (x_0\cdot a)=0,\, (x_0\cdot b)=0,$ and $\Xi\cap \mathbb{Z}^d\ne \emptyset$.
In addition, in this framework,
$\mathcal{P}f$ on $\mathpzc{T}^{\star}_{\alpha,\beta,a,b}$,
where $\beta\ge r$, uniquely determines $f$ on $\mathpzc{A}_{\alpha,\beta,a,b}$.
\label{thm2p2}
\end{theorem}

Note that $S_J$ is different from the origin $O$, that is $x_0$, 
in formula \eqref{defXS1} adopted for the plane
$\Xi=\Xi_{a,b}+x_0$
if $x_0\notin \mathbb{Z}^d$
in Theorem \ref{thm2p2}.

Proof of Theorem \ref{thm2p2} is similar to the proof of Theorem \ref{thm2p1} in view of the equalities that
$(\hat{z}_\star\cdot z)=0$ in \eqref{RTeqHz1snew},
$(\theta\cdot z)=0$ in \eqref{RTeqHz2new},
and the identity that
$\{x_0\}=\{x\in\Xi_{a,b}+x_0: x\cdot a=0, x\cdot b=0\}$, $\,x_0\in\mathbb{R}^d.$

Theorems \ref{thm2p1} and \ref{thm2p2} can be considered as discrete analogues of the classical slice-by-slice reconstruction for the continuous X-ray transform $\mathcal{P}$ in \eqref{XTeq1}.

In addition, Theorems \ref{thm2p1} and \ref{thm2p2}, for $d=2$, can be also considered as discrete analogues of the Cormack inversion \cite{Cormack}
for the continuous X-ray transform $\mathcal{P}$ in \eqref{XTeq1} for $d=2$.

In addition to Theorems \ref{thm2p1} and \ref{thm2p2}, where
$f$ is compactly supported, we also have the following conjectures.
\begin{conj}
Suppose that, for $\epsilon>0,$
\begin{equation}
f(x)=O(|x|^{-(1+\epsilon)}),\,\, x\in \mathbb{Z}^2,\,
\text{ as }|x|\to\infty.
\label{newassum}
\end{equation}
Then
$\mathcal{P}f$ on $\mathpzc{T}^\ast$ uniquely determines $f,$ where $\mathpzc{T}^\ast$ is
defined by \eqref{defTstar}
for $d=2$.
\label{conjec2p31}
\end{conj}
\begin{conj}
Suppose that 
\begin{equation}
f(x)=O(|x|^{-\infty}),\,\, x\in \mathbb{Z}^2,\,
\text{ as }|x|\to\infty.
\label{newassum}
\end{equation}
Then
$\mathcal{P}f$ on $\mathpzc{T}^\ast_{\alpha,\beta}$ uniquely determines $f$ on $\mathpzc{A}_{\alpha,\beta},$ where $\mathpzc{T}^\ast_{\alpha,\beta}$ and $\mathpzc{A}_{\alpha,\beta}$ are defined by \eqref{defCab},\eqref{defTstarab} for $d=2$,
and $\beta=\infty$.
\label{conjec2p3}
\end{conj}

Conjecture \ref{conjec2p31} and Conjecture \ref{conjec2p3} are formulated as discrete analogues of well known results in the continuous case; see, for example, \cite{Natterer}.

One can see that,
in Theorems \ref{thm2p1} and \ref{thm2p2},
the number of rational straight lines $\gamma$ used
for reconstructing $f$ 
on $B_r\cap\mathbb{Z}^d$
from $\mathcal{P}f(\gamma)$
is $N_r$ as in \eqref{defNrn}.
For Theorem \ref{thm2p1} and $d=2$, this reconstruction has a considerable similarity with the algorithmic reconstruction suggested in \cite{Kawazoe} (Section 4.3) for the case of function $f$ supported in a fixed square $[-r,r]^2$.
The point is that the latter reconstruction is also non-overdetermined although this important fact is not mentioned in \cite{Kawazoe}.
An essential advantage of our reconstructions in Theorems \ref{thm2p1} and \ref{thm2p2} is that the family of rational straight lines required for these reconstructions is given in advance by explicit analytical formulas in contrast with inductive definition of an appropriate family of straight lines in \cite{Kawazoe}. 
Moreover, our basic family of straight lines $\mathpzc{T}^\ast$ is defined for the whole space $\mathbb{Z}^d$
and is independent of $r.$
In particular, $\mathpzc{T}^\ast$ is used in Conjecture \ref{conjec2p31} and $\mathpzc{T}^\ast_{\alpha,\infty}$ is used in Conjecture \ref{conjec2p3} for $d=2$, where the parameter $r$ is absent since $f$ is not compactly supported.

\subsection{Weighted case}
\label{secwt2p1}
In addition to the X-ray transform $\mathcal{P}$ in \eqref{XTeq1}, and its discrete analogue in \eqref{XTeq1disc}, we also consider their weighted versions defined by the formulas:
\begin{equation}
\mathcal{P}_W f(\gamma)=\int_{\gamma}W(y,\hat{\gamma})\, f(y)dy,\quad \gamma\in T,
\label{XTeq1wt}
\end{equation}
for the continuous case;
\begin{equation}
\mathcal{P}_W f(\gamma)=\sum_{y\in\gamma\cap\mathbb{Z}^d} W(y,\hat{\gamma})\,f(y),\quad \gamma\in {T}',
\label{XTeq1discwt}
\end{equation}
for the discrete case.
Here, $W$ is a weight function, $T$ and $T'$ are defined by \eqref{XTeq2} and \eqref{defTTz}, and $\hat{\gamma}$ is the direction of $\gamma$.

We assume that 
\begin{equation}
\text{$W$ has no zeros}.
\label{assumW}
\end{equation}

Under this assumption,
the results in subsections \ref{sec2p1} and \ref{secrat} remain valid after minor modifications.
In place of formula \eqref{fdef1}, we have
\begin{equation}
f(x)=W(x,\theta)^{-1}\mathcal{P}_Wf(\gamma_{x,\theta}), \quad x\in\mathbb{Z}^d, \quad\theta\in \mathbb{S}^{d-1}\setminus\Omega.
\label{fdef122}
\end{equation}
Moreover, formulas \eqref{fdeff2}, \eqref{fdeff2n} admit similar generalizations.

In place of formulas \eqref{defrecon1}, \eqref{defrecon2}, we have
\begin{equation}
f(z)=W(z,\hat{\gamma}_z)^{-1}\mathcal{P}_Wf(\gamma_z), \quad z\in S_1,
\label{defrecon122}
\end{equation}
\begin{equation}
f(z)=W(z,\hat{\gamma}_z)^{-1}\left(\mathcal{P}_Wf(\gamma_z)-\sum_{\zeta\in \gamma_z\cap (\cup_{k=1}^i S_k)}W(\zeta,\hat{\gamma}_z)f(\zeta)\right), \quad z\in S_{i+1},
\label{defrecon222}
\end{equation}
for $i=1, 2, \dotsc, J-1,$
where
$\gamma_z$ is defined by \eqref{RTeqHz1}-\eqref{RTeqHz2}.

Taking into account formulas \eqref{defrecon122}, \eqref{defrecon222},
Theorems \ref{thm2p1}, \ref{thm2p2} admit straightforward generalizations for the case of weighted discrete X-ray transform $\mathcal{P}_W$ under assumption \eqref{assumW}.

In this respect, some results to the discrete weighted X-ray transform $\mathcal{P}_W$ in \eqref{XTeq1discwt} are drastically simpler than those for the continuous case; see, for example, \cite{Boman}, \cite{Goncharov} and references therein.

\section{Relations with the continuous case}
\label{sec3}
Let $\mathcal{P}^{con}$ denote $\mathcal{P}$ defined by \eqref{XTeq1}
and
$\mathcal{P}^{dis}$ denote $\mathcal{P}$ defined by \eqref{XTeq1disc}.
To present relations between $\mathcal{P}^{con}$ and $\mathcal{P}^{dis}$, we use
\begin{equation}
X_{\theta}=\{x\in\mathbb{R}^d: x\cdot\theta=0\}, \quad\theta\in\mathbb{S}^{d-1},
\end{equation}
the orthogonal projection $\pi_\theta$ of $\mathbb{R}^d$ to $X_{\theta}$, and 
\begin{equation}
X_{\theta}^{dis}=\pi_{\theta}\mathbb{Z}^d.
\end{equation}
Note that
\begin{equation}
X_{\theta}=\pi_{\theta}\mathbb{R}^d.
\end{equation}

\begin{remark}
The set $X_\theta^{dis}$ 
has the property that dist$(x,y)>\varepsilon>0$ for $x,y\in X_\theta^{dis},\, x\ne y$, for fixed $\theta\in\Omega$, where $\Omega$ is defined by \eqref{defOmg}.
\label{rem3p1}
\end{remark}
For completeness of presentation, this remark is proved in subsection \ref{pfrem3p1}.

On the other hand, the structure of $X_\theta^{dis}$ is more complicated for $\theta\in\mathbb{S}^{d-1}\setminus\Omega$. For example,
if $\theta_1, \dotsc, \theta_d$ are $\mathbb{Q}$-linearly independent,
then
$X_\theta^{dis}$ is everywhere dense in $X_\theta$
in view of the Kronecker-Weyl theorem.
See, for example, \cite{Bailleul} in connection with this theorem.

Note also that definitions \eqref{XTeq1} and \eqref{XTeq1disc} can be re-written as
\begin{equation}
\mathcal{P}^{con}f(\theta,x)=\int_{\mathbb{R}} f(x+s\theta)ds,\, x\in X_\theta, \,\theta\in \mathbb{S}^{d-1},
\label{XTeqnew}
\end{equation}
and
\begin{equation}
\mathcal{P}^{dis}f(\theta,x)=\sum_{s\in\mathbb{R}, \,x+s\theta\in\mathbb{Z}^d} f(x+s\theta),\, x\in X_\theta^{dis}, \,\theta\in \mathbb{S}^{d-1},
\label{XTeq1discnew}
\end{equation}
respectively.

\subsection{Case of delta functions}
\label{sec3p1}
In this subsection, we present relations
between 
$\mathcal{P}^{con}$ and $\mathcal{P}^{dis}$
for the case when $f$ in \eqref{XTeqnew} 
is a sum of Dirac-delta functions supported on $\mathbb{Z}^d\subset\mathbb{R}^d.$
\begin{theorem}
Let $f^{dis}$ be a function on $\mathbb{Z}^d$,
and
\begin{equation}
f(x)=\sum_{z\in\mathbb{Z}^d}f^{dis}(z)\delta(x-z), \quad x\in\mathbb{R}^d,
\label{fexpn}
\end{equation}
where $\delta$ denotes the Dirac-delta function.
Suppose that $f^{dis}$ satisfies condition \eqref{suppfBr}.
Then the following formula holds:
\begin{equation}
(\mathcal{P}^{con}f)_\theta(x)
=\sum_{\zeta\in X_{\theta}^{dis}}(\mathcal{P}^{dis}f^{dis})_{\theta}(\zeta)\delta(x-\zeta), \quad x\in X_\theta,\quad \theta\in\mathbb{S}^{d-1},
\label{themformula}
\end{equation}
where
\begin{equation}
(\mathcal{P}^{con}f)_\theta(x)=\mathcal{P}^{con}f(\theta,x),
\quad
(\mathcal{P}^{dis}f^{dis})_{\theta}(y)=\mathcal{P}^{dis}f^{dis}(\theta,y).
\label{propeqn}
\end{equation}
\label{pro3p1}
\end{theorem}

\subsection*{Proof of Theorem \ref{pro3p1}.}
For the continuous case, in view of \eqref{XTeqnew}, we have that
\begin{equation}
(\mathcal{P}^{con}\delta_{y})_\theta(x)
=\delta_{\pi_\theta y}(x),\quad 
x\in X_\theta,\, y\in\mathbb{R}^d,
\label{XTeqx1}
\end{equation}
where
\begin{equation}
 \delta_{\zeta}(x)=\delta(x-\zeta), \quad x, \zeta\in\mathbb{R}^d.
\end{equation}

For the discrete case, in view of \eqref{XTeq1discnew}, we have that
\begin{equation}
(\mathcal{P}^{dis}\tilde{\delta}_{z})_\theta(x)
=K_{\pi_\theta z}(x),\quad 
x\in X^{dis}_\theta,\, z\in\mathbb{Z}^d,
\label{XTeqx2}
\end{equation}
where 
\begin{equation}
\tilde{\delta}_z(x)=1\text{ if }x=z\text{ and }\tilde{\delta}_z(x)=0\text{ if }x\ne z, \,x, z\in\mathbb{Z}^d,
\end{equation}
\begin{align}
K_\zeta(x)=K(x-\zeta), \, x, \zeta\in\mathbb{R}^d, \quad 
K(x)=1\text{ if }x=0\text{ and }K(x)=0\text{ if }x\ne 0,
\end{align}
that is, $\tilde{\delta}$ and $K$ are analogues of Kronecker symbol.

In view of 
\eqref{XTeqx1},
\eqref{XTeqx2},
 we have that
\begin{align}
(\mathcal{P}^{con}\delta_y)_\theta(x)
&=\delta(x-\pi_\theta y)=\sum_{\zeta\in X_{\theta}^{dis}}K_{\pi_\theta y}(\zeta)\delta(x-\zeta)\label{thmderv}\\
&=\sum_{\zeta\in X_{\theta}^{dis}}(\mathcal{P}^{dis}\tilde{\delta}_y)_{\theta}(\zeta)\delta(x-\zeta),
\,\,x\in X_\theta,\, y\in\mathbb{R}^d.\notag
\end{align}

Formula \eqref{themformula} follows from \eqref{fexpn}, the linearity of the transforms $\mathcal{P}^{con}$ and $\mathcal{P}^{dis}$, and from Eq. \eqref{thmderv}.

Theorem \ref{pro3p1} is proved.

 In addition, in connection with our assumption \eqref{suppfBr}, we consider 
 \begin{equation}
 X_{\theta,r}^{dis}=\pi_{\theta}(B_r\cap\mathbb{Z}^d), \quad\theta\in\mathbb{S}^{d-1}, \, r>0.
 \end{equation} 
One can see that 
 \begin{equation}
 \#X_{\theta,r}^{dis}\le N_r, \quad\theta\in\mathbb{S}^{d-1}, \, r>0,
 \label{Xrt}
 \end{equation}
 \begin{equation}
 \#X_{\theta,r}^{dis}=N_r, \quad\theta\in\mathbb{S}^{d-1}\setminus\Omega, \, r>0,
 \label{Xrt2}
 \end{equation}
 where $N_r$ is defined by \eqref{defNrn}.

Let
\begin{equation}
B_{\theta,y,\epsilon}=\{x\in X_\theta: \text{dist }(x,y)<\epsilon\},\, y\in X_\theta, \epsilon>0,
\label{B1def}
\end{equation}
and
\begin{equation}
X_{\theta,r,\epsilon}=\{x\in X_{\theta}: \text{dist }(x,X^{dis}_{\theta,r})<\epsilon\},\,\epsilon>0.
\label{B12def}
\end{equation}

\begin{prop}
Let the assumptions of Theorem \ref{pro3p1} be satisfied.
Then the projection $(\mathcal{P}^{con}f)_\theta$ 
on $X_{\theta,r,\epsilon}$
uniquely determines $(\mathcal{P}^{dis}f^{dis})_\theta$, for any fixed $\theta\in\mathbb{S}^{d-1}$ and any $\epsilon>0$.
More precisely, for
each
$y\in X^{dis}_{\theta,r}$,
the projection $(\mathcal{P}^{con}f)_\theta$ 
on $B_{\theta,y,\epsilon}$
uniquely determines $(\mathcal{P}^{dis}f^{dis})_\theta(y)$, for any fixed $\theta\in\mathbb{S}^{d-1}$ and any $\epsilon>0$.
\label{prop32}
\end{prop}
Proposition \ref{prop32} follows 
from formulas \eqref{themformula}, \eqref{Xrt},
and the property that
\begin{equation}
C\delta(x-y)\text{ on }
B_{\theta,y,\epsilon}
\text{ uniquely determines }C\text{ for fixed }y, \theta, \epsilon,
\label{diracC}
\end{equation}
where $\delta$ is the Dirac-delta function on $X_\theta.$

\begin{cor}
 Let the assumptions of Theorem \ref{pro3p1} be satisfied,
 and $\#X_{\theta,r}^{dis}=N_r$ for $\theta\in\mathbb{S}^{d-1}$.
Then the projection $(\mathcal{P}^{con}f)_\theta$ on $X_{\theta,r,\epsilon}$ uniquely determines $f$, for any $\epsilon>0$.
\label{cor3p4}
\end{cor}

Corollary \ref{cor3p4} follows from
Proposition \ref{prop32}
and formula \eqref{fdef1}.
The point is that formula \eqref{fdef1} is always valid if $\#(X^{dis}_{\theta,r})=N_r,$ even if $\theta\in\Omega.$

Let
\begin{equation}
\Gamma_{z,\epsilon}=\{\gamma\in T: \hat{\gamma}=\hat{\gamma}_z,\, \text{dist }(\gamma,\gamma_z)<\epsilon\}, \, z\in\mathbb{Z}^d,
\end{equation}
\begin{equation}
\mathpzc{T}^{\star}_{\epsilon}=\cup_{z\in \mathbb{Z}^d}\Gamma_{z,\epsilon},
\label{defTstarcon}
\end{equation}
\begin{equation}
\mathpzc{T}^{\star}_{\alpha,\beta,a,b,\epsilon}=\cup_{z\in\mathpzc{A}_{\alpha,\beta,a,b}}\Gamma_{z,\epsilon},\, 0\le\alpha<\beta,
\label{defTstarabcon}
\end{equation}
where $\epsilon>0,$
$\gamma_z$ is defined as in \eqref{RTeqHz1}, \eqref{RTeqHz1new},
$\hat{\gamma}$ denotes the direction of $\gamma$,
$\mathpzc{A}_{\alpha,\beta,a,b}$ is defined as
in \eqref{defCab1},
and $a,b$ are the same as in \eqref{defXi}.

\begin{cor}
 Let the assumptions of Theorem \ref{pro3p1} be satisfied.
 Then
 $\mathcal{P}^{con}f$
 on $\mathpzc{T}^{\star}_{\epsilon}$
 uniquely determines $f$, for any $\epsilon>0$.
In addition, in this framework,
$\mathcal{P}^{con}f$ on $\mathpzc{T}^{\star}_{\alpha,\beta,a,b,\epsilon}$,
where $\beta\ge r$ and $r$ is the number in \eqref{suppfBr}, uniquely determines $f(x)$ for $x\in \mathbb{R}^d, |x|\ge \alpha$.
\label{cor3p5}
\end{cor}

Corollary \ref{cor3p5} follows from
\eqref{diracC} and Theorems \ref{thm2p1}, \ref{thm2p2}.

\subsection{Regular case}
\label{sec3p2}
In this subsection, we present relations
between 
$\mathcal{P}^{con}$ and $\mathcal{P}^{dis}$
for the case when $f$ in \eqref{XTeqnew} 
is regular on $\mathbb{R}^d$, e.g., piecewise constant.
We start with
\begin{equation}
f(x)=\sum_{z\in\mathbb{Z}^d\cap B_r}f^{dis}(z)w_z\chi_{z}(x-z), \quad x\in\mathbb{R}^d,
\label{fexpnnew}
\end{equation}
where
$f^{dis}$ satisfies \eqref{suppfBr},
$\chi_z$ are the characteristic functions of $B_{\rho_z}$,
$\rho_z$ are positive, 
$w_z$ are non-zero complex numbers, 
and
\begin{equation}
(B_{\rho_{z_1}}+z_1)\cap (B_{\rho_{z_2}}+z_2)=\emptyset\text{ if }z_1\ne z_2,
\label{pro3p1assum}
\end{equation}
where $B_r$ is defined by \eqref{diskBr}.
Let
\begin{align}
T''=\{\gamma\in T': \gamma\cap (B_{\rho_z}+z)=z\text{ or }\gamma\cap (B_{\rho_z}+z)=\emptyset,\,\forall z\in\mathbb{Z}^d\cap B_r\},
\end{align}
where $T'$ is defined by \eqref{defTTz}.

\begin{prop}
Under assumptions \eqref{fexpnnew},
\eqref{pro3p1assum},
the following formula holds:
    \begin{equation}
 \mathcal{P}^{con}f(\gamma)=
\mathcal{P}^{dis}_Wf^{dis}(\gamma) 
=\mathcal{P}^{dis}(Wf^{dis})(\gamma)
\text{ for } \gamma\in T'',
\label{formla76}
\end{equation}
where $\mathcal{P}^{dis}_W$ is as in \eqref{XTeq1discwt},
\begin{align}
W(y,\hat{\gamma})=W(y),
\\
 W(y)=2\rho_y w_y\text{ if }y\in\gamma\cap\mathbb{Z}^d\cap B_r,\quad
 W(y)=1\text{ if }y\in\gamma\cap\mathbb{Z}^d, |y|>r,
\end{align}
and $\hat{\gamma}$ denotes the direction of $\gamma.$
\label{propnew3p6}
\end{prop}
Proposition \ref{propnew3p6} follows from the definitions of $\mathcal{P}^{con}$ and $\mathcal{P}^{dis}$
and is the simplest analogue of Theorem \ref{pro3p1} for the regular case.

However, for the regular case, some more interesting relations between continuous and discrete X-ray transforms are as follows.

Note that $\mathbb{R}^d$ can be presented as 
\begin{equation}
 \mathbb{R}^d=\cup_{z\in\mathbb{Z}^d}
 {\mathcal{U}}_z,
 \label{news1}
\end{equation}
 \begin{equation}
 \mathcal{U}_z=\mathcal{U}+z,\, \mathcal{U}=\{x=(x_1, \dotsc, x_d)\in\mathbb{R}^d: -1/2<x_j\le 1/2, \, j=1, \dotsc, d\}.
 \end{equation}

A regular function $f^{reg}$ supported in $B_r$ can be approximated by its piecewise constant version $f$ such that 
\begin{equation}
 f(x)=c(z), \,x\in\mathcal{U}_z, \, z\in \mathbb{Z}^d,
 \label{fpiece}
\end{equation}
where $c(z)=0$ if $z\in\mathbb{Z}^d\setminus B_r;$
e.g., one can assume that
$c(z)=f^{reg}(z), z\in \mathbb{Z}^d.$
In addition to $f$
in 
\eqref{fpiece},
we also consider
\begin{equation}
 f^{dis}(z)=c(z),\,z\in \mathbb{Z}^d.
 \label{fdisc}
\end{equation}
One can see that
\begin{equation}
f(x)=\sum_{\zeta\in \mathbb{Z}^d}f^{dis}(\zeta)\chi_{\zeta}(x), \, x\in\mathbb{R}^d,
\label{newformPcon2}
\end{equation}
where
$\chi_{\zeta}$ is the characteristic function of $\mathcal{U}_{\zeta}$.

We set
\begin{equation}
W(\zeta,\theta)
:=\mathcal{P}^{con}\chi_{\zeta}(\gamma)
=|\gamma\cap {\mathcal{U}_\zeta} |, \, \zeta\in\gamma, \hat{\gamma}=\theta, \, \zeta\in\mathbb{Z}^d, \theta\in\mathbb{S}^{d-1},\, \gamma\in T',
\label{defWnew}
\end{equation}
where 
$|\cdot|$ denotes the length
and
$T'$ is defined by \eqref{defTTz}.

\begin{prop}
Under assumptions \eqref{news1}-\eqref{fdisc},
the following formula holds:
    \begin{align}
\mathcal{P}^{dis}_Wf^{dis}(\gamma)
=\mathcal{P}^{con}f(\gamma)
-\mathcal{P}^{con}(\delta f)(\gamma), \quad \gamma\in T',
\label{recon1anew}
\end{align}
where 
\begin{equation}
\delta f(x,\gamma)
 =\sum_{\zeta\in \mathbb{Z}^d\setminus\gamma}f^{dis}(\zeta)\chi_{\zeta}(x), \, x\in\mathbb{R}^d,
\label{recon1anewaa}
\end{equation}
$\mathcal{P}^{dis}_W$ is defined according to \eqref{XTeq1discwt},
\eqref{defWnew},
$\chi_{\zeta}$ is the characteristic function of $\mathcal{U}_{\zeta}$,
and $T'$ is defined by \eqref{defTTz}.
\label{prop3p7}
\end{prop}

Proposition \ref{prop3p7} can be also considered as an analogue of Theorem \ref{pro3p1} for the regular case.
 
Proposition \ref{prop3p7} follows from \eqref{XTeq1discwt}, \eqref{news1}-
\eqref{defWnew}, and the formulas:
\begin{equation}
 \delta f(x,z)
 =f(x)-\sum_{\zeta\in\gamma\cap \mathbb{Z}^d}f^{dis}(\zeta)\chi_{\zeta}(x)
 =\sum_{\zeta\in \mathbb{Z}^d\setminus\gamma}f^{dis}(\zeta)\chi_{\zeta}(x),
\label{recon2a}
\end{equation}
\begin{equation}
 (\mathcal{P}^{con}(\sum_{\zeta\in\gamma\cap \mathbb{Z}^d}f^{dis}(\zeta)\chi_{\zeta})(\gamma)=
 \sum_{\zeta\in\gamma\cap \mathbb{Z}^d}f^{dis}(\zeta)W(\zeta,\hat{\gamma}).
\label{newformPcon1}
\end{equation}

Next, we consider the  problem of finding $f$ in \eqref{fpiece}
from $g=\mathcal{P}^{con}f$ on $\mathpzc{T}^\star$, where $\mathpzc{T}^\star$ is as in Theorems \ref{thm2p1}, \ref{thm2p2}.

Below we suggest 
a
reconstruction of $f$ from
$g$
on $\mathpzc{T}^\star$
via iterative use of discrete reconstructions in Theorems \ref{thm2p1}, \ref{thm2p2}
and
Proposition \ref{prop3p7}.

We present this construction next for $d=2$.

Proposition \ref{prop3p7}
implies the following Corollary.
\begin{cor}
Under assumptions \eqref{news1}-\eqref{fdisc},
for $d=2,$
the following formula holds for $f$ supported in $B_r$:
\begin{align}
\mathcal{P}^{dis}_Wf^{dis}(\gamma_z)
&=g(\gamma_z)-\sum_{\zeta\in \cup_{i<j} S_i\setminus\gamma_z}|\gamma_z\cap \mathcal{U}_\zeta|f^{dis}(\zeta)\notag\\
&-\sum_{\zeta\in \cup_{i\ge j} S_i\setminus\gamma_z}|\gamma_z\cap \mathcal{U}_\zeta|f^{dis}(\zeta), \quad z\in S_{j},
\label{defrecon3}
\end{align}
where $\gamma_z$ is defined by \eqref{RTeqHz1}-\eqref{RTeqHz2},
$S_j$ are defined by \eqref{defXS1}, \eqref{defXSj}, $j=1, 2, \dotsc, J.$
    \label{cor3p8}
\end{cor}

Our aforementioned iterative reconstruction of $f$ from 
$g$ on $\mathpzc{T}^\star$
is motivated by formula
\eqref{defrecon3}, 
and
the observation that 
\begin{align}
\sum_{\zeta\in \cup_{i\ge j} S_i\setminus\gamma_z}|\gamma_z\cap \mathcal{U}_\zeta|f^{dis}(\zeta)\approx0.
\label{defrecon3bb}
\end{align}
Here, we use that
if $z\in S_j, \, j=1, 2, \dotsc, J$, defined by \eqref{defXS1}, \eqref{defXSj},
$\zeta\in
\cup_{i\ge j}S_i$,
$\zeta\ne z,$
and $\gamma_z$ is defined by \eqref{RTeqHz1},
then
$\gamma_z\cap \mathcal{U}_\zeta$
is empty or, at least, typically, its length $|\gamma_z\cap \mathcal{U}_\zeta|$ is much smaller than $|\gamma_z\cap \mathcal{U}_z|$
even 
if
$\zeta$ is one of the nearest points to $z$ on $\cup_{i\ge j}S_i$.

Using \eqref{defrecon3}, \eqref{defrecon3bb}, we get
\begin{equation}
\mathcal{P}_W^{dis}f^{dis}(\gamma_z)\approx g(\gamma_z) ,\, z\in S_1,
\label{eq87ref}
\end{equation}
and, in general,
\begin{align}
\mathcal{P}^{dis}_Wf^{dis}(\gamma_z)
\approx
g(\gamma_z)-\sum_{\zeta\in \cup_{i<j} S_i\setminus\gamma_z}|\gamma_z\cap \mathcal{U}_\zeta|f^{dis}(\zeta), \quad z\in S_{j},\label{defrecon3a}
\end{align}
for $j=1, 2, \dotsc, J.$
In turn, formulas
\eqref{XTeq1discwt},
\eqref{defrecon122},
\eqref{defrecon222},
\eqref{defWnew}, 
\eqref{defrecon3a}
imply that
\begin{align}
w(z)f^{dis}(z)
\approx
g(\gamma_z), \quad z\in S_{1},\label{defrecon2itr4a}\\
w(z)f^{dis}(z)
\approx
g(\gamma_z)
-\sum_{\zeta\in \cup_{i<j} S_i}|\gamma_z\cap \mathcal{U}_\zeta|f^{dis}(\zeta), \quad z\in S_{j},
\label{defrecon2itr4}
\end{align}
for $j=1, 2, \dotsc, J,$
where
\begin{equation}
w(z)=W(z,\hat{\gamma}_z)=|\gamma_z\cap\mathcal{U}_z|.
\label{defwf}
\end{equation}

Here, $g(\gamma_z)=(\mathcal{P}^{con} f)(\gamma_z)$ for $f$ defined as in \eqref{fpiece}, \eqref{fdisc}, \eqref{newformPcon2}.

The point is that formulas \eqref{suppfBr}, \eqref{Brcov},
\eqref{defrecon2itr4a}, \eqref{defrecon2itr4}
yield
a 
recurrent layer-by-layer reconstruction procedure for $f^{dis}_1\approx f^{dis}$.

Improved approximations $f^{dis}_{i+1}$ can be constructed
as
follows.
First, due
to
\eqref{XTeq1discwt},
\eqref{recon1anew}, 
\eqref{recon1anewaa}, 
\eqref{defwf},
we have that
\begin{align}
(\mathcal{P}^{dis}(wf^{dis}))(\gamma_z)
=g(\gamma_z)-\sum_{\zeta\in \mathbb{Z}^d\setminus\gamma_z}f^{dis}(\zeta)|\gamma_z\cap \mathcal{U}_\zeta|, \, z\in \mathbb{Z}^d.
\label{defstep1}
\end{align}
Second, using formula \eqref{defstep1} and the approximation $f^{dis}\approx f^{dis}_i,$
we have that
\begin{equation}
(\mathcal{P}^{dis}(wf^{dis}))(\gamma_z)
\approx g^{dis}_i(\gamma_z):=
g(\gamma_z)-\sum_{\zeta\in \mathbb{Z}^d\setminus\gamma_z}f^{dis}_i(\zeta)|\gamma_z\cap \mathcal{U}_\zeta|, \, z\in \mathbb{Z}^d.
\end{equation}
Finally,
\begin{equation}
f^{dis}_{i+1}=w^{-1}(\mathcal{P}^{dis})^{-1}g^{dis}_i,
\end{equation}
where $(\mathcal{P}^{dis})^{-1}$ is realized
as in Theorem \ref{thm2p1}
via formulas
\eqref{defrecon1}, \eqref{defrecon2}.

It remains to recall
that $f$ and $f^{dis}$ are related by formula
\eqref{newformPcon2}.

\begin{remark}
Further 
theoretical and numerical
studies of applications 
of Theorem \ref{thm2p1}
and
Theorem \ref{thm2p2}
to reconstruction of $f$
from $\mathcal{P}^{con}f$
for the case of regular $f$ as in \eqref{fpiece}, \eqref{newformPcon2}
will be given in future works.
\end{remark}

\appendix
\section{Proofs of some estimates}
\label{secproof}

\subsection{Proof of estimate \eqref{defT}}
\label{secproof4p1}

Formula \eqref{defT} is equivalent to the 
following estimates:
\begin{equation}
C_1 r^{2d}< \# \mathpzc{T}_r^{\min}<C_2 r^{2d},
\label{refTr0min1}
\end{equation}
for some positive $C_1, C_2$
 depending on $d$.
The upper bound in
\eqref{refTr0min1},
with a non-optimal $C_2$, follows from the observation that $\# \mathpzc{T}_r^{\min}<N_r^2$, where $N_r$ is defined in \eqref{defNrn}.
The lower bound in \eqref{refTr0min1} can be proved as follows.

Let 
\begin{equation}
\mathpzc{T}_{r,0}^{\min}:=\{\gamma\in \mathpzc{T}_r^{\min}: 0\in\gamma\}.
\end{equation}
We have that
\begin{equation}
c_1 r^{d}< \# \mathpzc{T}_{r,0}^{\min}<c_2 r^{d},
\label{refTr0min}
\end{equation}
for some positive $c_1, c_2$ depending on $d$.
The estimates \eqref{refTr0min} follow from the asymptotic formula
\begin{equation}
\#F_n \sim \frac{n^{d}}{d\,\zeta(d)},
\,\text{ as }n\to\infty,
\label{refTr0min23}
\end{equation}
where 
\begin{equation}\label{fareydef}
	F_n:=\bigg\{ \frac{\boldsymbol{p}}{q} \in[0,1)^{d-1} : (\boldsymbol{p},q)\in{\widehat{\mathbb Z}}^{d}, \; 0<q\leq n \bigg\},
\end{equation}
\begin{equation}\label{fareydef1}
{\widehat{\mathbb Z}}^{d}:=\{ {{\text{\boldmath$m$}}}\in{\mathbb Z}^{d}\setminus\{0\} : \gcd({{{\text{\boldmath$m$}}}})=1 \},
\end{equation}
and $\zeta$ is the Riemann zeta function;
see, for example, \cite{Marklof} and references therein.
The set $F_n$ is known as the Farey points of level $n.$
In particular, for $d=2$, 
$\#{F_n}\sim \frac{3n^2}{\pi^2}$ 
and
$F_n$
is the sequence of completely reduced fractions $z_2/z_1$,
where $z_1, z_2\in\mathbb{Z}, \,z_1\ne0,$ 
$0\le z_2< z_1\le n$.

To get \eqref{refTr0min},
one can use \eqref{refTr0min23},
with $n
=\lceil r\rceil$
and
$n=\lfloor r/\sqrt{d}\rfloor,$ respectively,
for the upper and lower bounds. In a way similar to \eqref{refTr0min1}, the upper bound in \eqref{refTr0min} also follows from the observation that $\#\mathpzc{T}_{r,0}^{\min}<N_r$.

We introduce
\begin{equation}
\mathpzc{T}_{r,z}^{\min}:=\{\gamma\in \mathpzc{T}: z\in\gamma, \#(\gamma\cap B_{r,z}\cap\mathbb{Z}^d)\ge2\},
\end{equation}
where $\mathpzc{T}$ is defined by \eqref{defTset},
\begin{equation}
B_{r,z}:=\{x\in\mathbb{R}^d: |x-z|\le r\},\,z\in\mathbb{Z}^d,\, r>0.
\end{equation}

Let
\begin{equation}
B_r\cap\mathbb{Z}^d=\cup_{j=0}^{N_r-1} z_j,\, z_j\in\mathbb{Z}^d,\, |z_j|\le |z_{j+1}|.
\end{equation}
Note that $z_0=0\in\mathbb{Z}^d.$

One can see that:
\begin{equation}
\# (\mathpzc{T}_{r,z_1}^{\min}\setminus \mathpzc{T}_{r,0}^{\min})=\#\mathpzc{T}_{r,z_1}^{\min}-1;
\label{newwq1a}
\end{equation}
\begin{equation}
\# (\mathpzc{T}_{r,z_2}^{\min}\setminus(\mathpzc{T}_{r,z_1}^{\min}\cup\mathpzc{T}_{r,0}^{\min}))\ge\#\mathpzc{T}_{r,z_2}^{\min}-2;
\end{equation}
$\# (\mathpzc{T}_{r,z_2}^{\min}\setminus(\mathpzc{T}_{r,z_1}^{\min}\cup\mathpzc{T}_{r,0}^{\min}))=\#\mathpzc{T}_{r,z_2}^{\min}-2$
if 
$z_2\in B_{r,z_1}$;
and, in general,
\begin{equation}
\# (\mathpzc{T}_{r,z_j}^{\min}\setminus \cup_{i=0}^{j-1}\mathpzc{T}_{r,z_i}^{\min})\ge\#\mathpzc{T}_{r,z_j}^{\min}-j;
\label{newwq2a}
\end{equation}
$\# (\mathpzc{T}_{r,z_j}^{\min}\setminus \cup_{i=0}^{j-1}\mathpzc{T}_{r,z_i}^{\min})=\#\mathpzc{T}_{r,z_j}^{\min}-j$
if $z_i\in \cap_{k=1}^{i-1}B_{r,z_{k}}$, $i=1, 2, \dotsc, j\le N_r-1$.

Also, one can see that:
\begin{equation}
\cup_{j=0}^{N_{r/2}-1} \mathpzc{T}_{r/2,z_j}^{\min}\subset\mathpzc{T}_{r}^{\min};
\label{newwq1}
\end{equation}
$|z_{j_1}-z_{j_2}|\le r/2, \,\text{i.e. }\,z_{j_2}\in B_{r/2,z_{j_1}}, \, j_1, j_2\le N_{r/4}-1.$
Using formula \eqref{newwq1}
and
formulas \eqref{newwq1a}, \eqref{newwq2a},
with $r$ replaced by $r/2$,
we get
\begin{equation}
\#\mathpzc{T}_{r}^{\min}>\sum_{j=0}^{N_{r/2}-1}(\#\mathpzc{T}_{r/2,z_j}^{\min}-j)=\frac{1}{2}N_{r/2}(1+N_{r/2}).
\label{newwq1ab}
\end{equation}
Formula \eqref{newwq1ab} implies the lower bound in \eqref{refTr0min1}, at least, with non-optimal $C_1$.

\subsection{Proof of Remark \ref{rem3p1}}
\label{pfrem3p1}
Let $\theta\in \Omega$ and $\theta=\zeta/|\zeta|$, where $\zeta\in\mathbb{Z}^d\setminus\{0\}$ is such that $|\zeta|$ is minimal for fixed $\theta$.
Then 
\begin{equation}
 \pi_\theta z=z-(z\cdot\zeta) \frac{\zeta}{|\zeta|^2},\,\, z\in\mathbb{Z}^d.
\end{equation}
It is sufficient to show that 
\begin{equation}
|\pi_\theta z|\ge\varepsilon(\zeta)=|\zeta|^{-1},
\label{pf3p1eq1}
\end{equation}
$\text{ for }z\in\mathbb{Z}^d,\, z\ne s\zeta, \, s\in\mathbb{R}$.

We have that
\begin{equation}
|\zeta|^2|\pi_\theta z|^2=||\zeta|^2 z-(z\cdot\zeta) \zeta|^2
=|\zeta|^4 |z|^2-|\zeta|^2(z\cdot\zeta)^2\in\mathbb{Z}\cap[0,\infty),
\label{pf3p1eq2}
\end{equation}
\begin{equation}
|\zeta|^2|\pi_\theta z|^2=|\zeta|^4 |z|^2-|\zeta|^2(z\cdot\zeta)^2>0, \text{ for }z\ne s\zeta, \, s\in\mathbb{R}.
\label{pf3p1eq3}
\end{equation}

Formula \eqref{pf3p1eq1} follows from \eqref{pf3p1eq2}, \eqref{pf3p1eq3}.

This completes the proof of Remark \ref{rem3p1}.

\section*{Acknowledgments}
A part of this work was completed during the stay (in Nov-Dec 2025) of the second author at the Institut des Hautes Études Scientifiques (IHES). The authors would like to thank  IHES and \'{E}cole polytechnique for support during the visit.

\end{document}